\title{$\mathbb{Z}_2^2$-actions on Horikawa surfaces.}

\author{Vicente Lorenzo}
\date{}
\RequirePackage{fix-cm}
\RequirePackage{amsmath}

\documentclass{article}
\usepackage{graphicx}
\usepackage[T1]{fontenc}
\usepackage[utf8]{inputenc}
\usepackage{lmodern}
\usepackage[english]{babel}
\usepackage{amssymb}
\usepackage{amsthm}
\usepackage[all]{xy}
\usepackage{enumerate}
\usepackage[none]{hyphenat}
\usepackage{etoolbox}
\usepackage{microtype}
\usepackage{lipsum}

\newtheorem{theorem}{Theorem}

\newtheorem{proposition}{Proposition}

\theoremstyle{remark}\newtheorem{remark}{Remark}
\newenvironment{acknowledgements}{\textit{Acknowledgements.}}{}

\newcommand\blfootnote[1]{%
  \begingroup
  \renewcommand\thefootnote{}\footnote{#1}%
  \addtocounter{footnote}{-1}%
  \endgroup
}

\begin{document}

\maketitle

\begin{abstract}
Minimal algebraic surfaces of general type $X$ such that $K^2_X=2\chi(\mathcal{O}_X)-6$ or 
$K^2_X=2\chi(\mathcal{O}_X)-5$ are called Horikawa surfaces. In this note we study $\mathbb{Z}^2_2$-actions
on Horikawa surfaces.
The main result is that all the connected components of Gieseker's moduli space of canonical
models of surfaces of general type with invariants satisfying these relations
contain surfaces with $\mathbb{Z}^2_2$-actions.
\end{abstract}


\section{Introduction.}
Let $X$ be a minimal surface of general type over the complex numbers, which will be the ground field throughout the paper. If we denote by $K^2_X$ the self intersection
of its canonical class and by $\chi(\mathcal{O}_X)$ its holomorphic Euler characteristic, it is well known
(cf. \cite{Barth}) 
that the following inequalities are satisfied
\begin{equation}\label{GeneralType}
 \chi(\mathcal{O}_X)\geq 1,\quad K_X^2\geq 1, \quad 2\chi(\mathcal{O}_X)-6\leq K_X^2\leq 9\chi(\mathcal{O}_X).
\end{equation}
A pair of integers $(\chi, K^2)$ is said to be an admissible pair if it satisfies 
the inequalities (\ref{GeneralType}).
 The geographical question, i.e.
 whether there exists a
minimal algebraic surface of general type $X$ with $(\chi(\mathcal{O}_X), K^2_X)=(\chi,K^2)$ 
for any admissible pair $(\chi, K^2)$, has received much attention
{{\cite[Section VII.8]{Barth}}}. Several authors have also studied the geography of surfaces with special features 
like genus $2$ fibrations (see \cite{Per2}), simply-connectedness (see \cite{Per1}), $2$-divisibility of the
canonical class (see \cite{PerssonGang}), global $1$-forms (see \cite{Mendespardini1}), etc.
In this note we are interested in the geography of surfaces with $\mathbb{Z}_2^2$-actions. 

Let us denote by $\mathfrak{M}_{\chi,K^2}$ Gieseker’s moduli space of canonical
models of surfaces of general type
with fixed self intersection of the canonical class $K^2$ and fixed holomorphic Euler characteristic $\chi$. Then,
for any admissible pair $(\chi, K^2)$,
the moduli space $\mathfrak{M}_{\chi,K^2}$
has a finite number of connected components which parametrize the deformation equivalence classes 
of  canonical
models of surfaces of general type (cf. \cite{CataDef}). 
The main result is the following:
\begin{theorem}\label{MainResult}
Let $(\chi,K^2)$ be a pair of integers such that $K^2=2\chi-6$ and $\chi\geq 4$ or 
$K^2=2\chi-5$ and $\chi\geq 3$. 
Then it is possible to find a surface with a $\mathbb{Z}_2^2$-action 
in any of the connected components of $\mathfrak{M}_{\chi,K^2}$.
\end{theorem}

\begin{remark}
It follows from (\ref{GeneralType}) that if we want $(\chi, K^2)$ to be the invariants of a minimal 
surface of general type and $K^2=2\chi-6$ (resp. $K^2=2\chi-5$),
then the inequality $\chi\geq 4$ (resp. $\chi\geq 3$) is not a restriction.
\end{remark}

The note is structured as follows. Section \ref{Horikawa surfaces.} is devoted to present some properties 
of minimal surfaces of general type $X$ with $K_X^2=2\chi(\mathcal{O}_X)-6$ and to study
their deformation equivalence classes.
In Section \ref{Horikawa surfaces II.} 
we do the same for minimal surfaces of general type $X$ with $K_X^2=2\chi(\mathcal{O}_X)-5$.
 In Section \ref{Z22 covers.}
how to construct $\mathbb{Z}_2^2$-covers and to obtain information about them is explained. 
Sections \ref{Constructions.} and \ref{ConstructionsSecondCase.} 
include the examples that prove Theorem \ref{MainResult}.
Finally, Section \ref{Genus2Fibrations} gives an insight into these examples. 

\section{Horikawa surfaces on the line $K^2=2\chi-6$.} \label{Horikawa surfaces.} 
If $X$ is an algebraic surface then $\chi(\mathcal{O}_X)=1-q(X)+p_g(X)$ where
$p_g(X)=h^0(K_X)$ is the geometric genus of $X$ and $q(X)=h^0(\Omega_X^1)$ is its irregularity.
If $X$ is minimal and of general type we know by {\cite[Proposition 2.3.2]{Mendespardini1}} that 
$q(X)>0$ implies $K^2_X\geq 2\chi(\mathcal{O}_X)$. On the other hand, if $K_X^2=2p_g(X)-4$
then $q(X)=0$ by {\cite[Theorem 10]{Bomb}}.

In \cite{Hor1} Horikawa studied minimal surfaces of general type $X$ such that $K_X^2=2p_g(X)-4$. We know by 
the previous comments that this is equivalent to $K_X^2=2\chi(\mathcal{O}_X)-6$.
The following theorems are some of the results Horikawa proved in \cite{Hor1}.

\begin{theorem}
[{{\cite[Lemma 1.1 and Lemma 1.2]{Hor1}}}]
Let $X$ be a minimal algebraic surface with $K_X^2=2\chi(\mathcal{O}_X)-6$
and $\chi(\mathcal{O}_X) \geq 4$. Then the canonical system $|K_X|$ has no base point. Moreover,
the canonical map $\varphi_{K_X}:X\rightarrow \mathbb{P}^{p_g(X)-1}$ is a morphism 
of degree 2 onto a surface that is isomorphic to either $\mathbb{P}^2$, a Hirzebruch surface $\mathbb{F}_e$ or
a cone over a rational curve of degree $e$ in $\mathbb{P}^e$.
\end{theorem}

\begin{remark}
 It is worth noticing that $\mathbb{P}^2$ can only be obtained if $\chi\in\{4,7\}$ and 
 a cone over a rational curve of degree $e$ in $\mathbb{P}^e$ can only be obtained if $\chi\in\{5,6,7\}$.
\end{remark}

\begin{remark}\label{PreGenus2Result}
Suppose that the canonical image of a minimal surface of general type $X$ with
$K_X^2=2\chi(\mathcal{O}_X)-6$ is a Hirzebruch surface $\mathbb{F}_e$. 
We denote by $\Delta_0$ the negative section of $\mathbb{F}_e$ and by $F$ the fiber with
respect to the projection $\mathbb{F}_e\rightarrow \mathbb{P}^1$. 
Then we can realize $X$ as a double cover  $X\rightarrow \mathbb{F}_e$ branched along a divisor of 
type $6\Delta_0+mF$ for some positive integer $m$. By the Riemann-Hurwitz formula a general fiber $F$
pulls back to a curve of genus $2$ in $X$.  It follows that the composition
$X\rightarrow \mathbb{F}_e\rightarrow \mathbb{P}^1$ is a genus $2$ fibration.

Suppose now that the canonical image of $X$ is a cone $\Sigma$ over a rational
curve of degree $e$ in $\mathbb{P}^e$. Then we can decompose the canonical map 
$X\rightarrow \Sigma\subset \mathbb{P}^{e+1}$ of $X$ as
$X\rightarrow \mathbb{F}_{e}\rightarrow \Sigma\subset \mathbb{P}^{e+1}$
where $X\rightarrow \mathbb{F}_{e}$ is a double cover and
$\mathbb{F}_{e}\rightarrow \Sigma\subset \mathbb{P}^{e+1}$ contracts the 
negative section of $\mathbb{F}_{e}$ and is an isomorphism elsewhere. 
Moreover the branch locus of $X\rightarrow \mathbb{F}_{e}$ is again of type 
$6\Delta_0+mF$ so we get a genus $2$ fibration on $X$ also in this case.
\end{remark}

We can summarize Remark \ref{PreGenus2Result} as follows.

\begin{proposition}[{{\cite[Corollary 1.7]{Hor1}}}]\label{Gen26}
Let $X$ be a minimal surface of general type with $K_X^2=2\chi(\mathcal{O}_X)-6$. If the canonical image of $X$
is not $\mathbb{P}^2$ then it admits a genus $2$ fibration $X\rightarrow \mathbb{P}^1$.
\end{proposition}

\begin{remark}\label{Case42}
Let $X$ be a smooth minimal surface with $\chi(\mathcal{O}_X)=4$ and $K^2_X=2$. Then the canonical image of
$X$ is $\mathbb{P}^2$ and we cannot apply Proposition \ref{Gen26}. Nevertheless, in Section \ref{Constructions.}
we are going to construct a surface with these invariants and a genus $2$ fibration. Now, if $X$ has a genus $2$ 
fibration then $K_X$ cannot be ample. Indeed, let us denote by $F$ a general fiber of this genus $2$ fibration.
The exact sequence
\begin{equation*}
 0\rightarrow \mathcal{O}_X(K_X-F)\rightarrow \mathcal{O}_{X}(K_X)\rightarrow \mathcal{O}_{F}(K_F)\rightarrow 0
\end{equation*}
induces an exact sequence
\begin{equation*}
 0\rightarrow H^0(X,K_X-F)\rightarrow H^0(X,K_X)
 \rightarrow H^0(F,K_F)\rightarrow \cdots
\end{equation*}
Since $h^0(X,K_X)=3$ and $h^0(F,K_F)=2$ we conclude that $H^0(X,K_X-F)\neq 0$. Let us consider $D\in |K_X-F|$.
Then $K_XD=0$ and $D^2=-2$. It follows that $D$ contains at least one $(-2)$-curve.

\end{remark}

\begin{theorem}[{{\cite[Theorem 3.3, Theorem 4.1 and Theorem 7.1]{Hor1}}}]
\label{DefClass}
Let $(\chi,K^2)$ be a pair of integers such that $K^2=2\chi-6$ and $\chi\geq 4$.
If $K^2\notin 8\cdot\mathbb{Z}$ then minimal algebraic surfaces $X$ such
that $(\chi(\mathcal{O}_X), K^2_X)=(\chi,K^2)$ have one and the same deformation type.
If $K^2\in 8\cdot\mathbb{Z}$ then minimal algebraic surfaces $X$ such 
that $(\chi(\mathcal{O}_X), K^2_X)=(\chi,K^2)$ have two deformation classes.
The image of the canonical map of a surface in the first class is $\mathbb{F}_e$
for some $e\in\{0,2,\ldots,\frac{1}{4}K^2\}$. The image of the canonical map of a 
surface in the second class is $\mathbb{F}_{\frac{1}{4}K^2+2}$ if $K^2>8$ and $\mathbb{P}^2$
or a cone over a rational curve of degree $4$ in $\mathbb{P}^4$ if $K^2=8$.
\end{theorem}

\section{Horikawa surfaces on the line $K^2=2\chi-5$.} \label{Horikawa surfaces II.}
In \cite{Hor2} Horikawa studied minimal surfaces of general type $X$ such that $K_X^2=2p_g(X)-3$. 
Arguing like we did in Section \ref{Horikawa surfaces.}, it follows from {\cite[Proposition 2.3.2]{Mendespardini1}} and {\cite[Theorem 10]{Bomb}}
that this is equivalent to $K_X^2=2\chi(\mathcal{O}_X)-5$.
The following theorems are some of the results Horikawa proved in \cite{Hor2}.

\begin{theorem}
[{{\cite[Section 1]{Hor2}}}]\label{CanSys5}
Let $X$ be a minimal algebraic surface with $K_X^2=2\chi(\mathcal{O}_X)-5$
and $\chi(\mathcal{O}_X) \geq 6$. Then the canonical system $|K_X|$ has a unique base point $b$ and
 the canonical map $\varphi_{K_X}:X\dashrightarrow \mathbb{P}^{p_g(X)-1}$
is a rational map whose image is isomorphic to either a Hirzebruch surface $\mathbb{F}_e$ or
a cone over a rational curve of degree $e$ in $\mathbb{P}^e$. Moreover, let us denote by 
$\varepsilon:\widetilde{X}\rightarrow X$ the blow-up of $X$ at $b$. Then $\widetilde{X}$ is the minimal
resolution of singularities of the normal double cover $X'$ of one of the following surfaces.
\begin{enumerate}
\item[a)] A blow-up of the Hirzebruch surface $\mathbb{F}_e$ at two points of a fiber that may be infinitely near.
\item[b)] A blow-up of the Hirzebruch surface $\mathbb{F}_1$ at a point on the negative section.
\item[c)] The Hirzebruch surface $\mathbb{F}_2$.
\end{enumerate}
\end{theorem}

\begin{remark}
It is worth noticing that assuming $\chi\geq 6$ the canonical image can be isomorphic to a cone over a
rational curve only when $(\chi, K^2)=(6,7)$. In addition, case $b)$ of Theorem \ref{CanSys5} can only happen
when $(\chi, K^2)=(6,7)$ and
case $c)$ of Theorem \ref{CanSys5} can only happen when $(\chi, K^2)=(7,9)$.
\end{remark}

Similarly to Remark \ref{PreGenus2Result} and using the notation of Theorem
 \ref{CanSys5} it can be proved that:

\begin{proposition}[{{\cite[Corollary 1.4]{Hor2}}}]\label{Gen25}
Let $X$ be a minimal surface of general type with $K_X^2=2\chi(\mathcal{O}_X)-5$.
If $X'$ is a double cover of the blow-up of $\mathbb{F}_e$ at two points of a fiber, 
then $X$ admits a genus $2$ fibration $X\rightarrow \mathbb{P}^1$.
\end{proposition}

\begin{theorem}[{{\cite[Introduction]{Hor2}}}]
\label{DefClass2}
Let $(\chi,K^2)$ be a pair of integers such that $K^2=2\chi-5$ and $\chi\geq 3$.
If $K^2+1\notin 8\cdot\mathbb{Z}$ and $(\chi,K^2)\neq (7,9)$ then minimal algebraic
surfaces $X$ such that $(\chi(\mathcal{O}_X), K^2_X)=(\chi,K^2)$ have one and the same deformation type.
If $K^2+1\in 8\cdot\mathbb{Z}$ or $(\chi,K^2)=(7,9)$ then minimal algebraic surfaces $X$ such that 
$(\chi(\mathcal{O}_X), K^2_X)=(\chi,K^2)$ have two deformation classes.
\end{theorem}

\begin{remark}
Let $(\chi,K^2)$ be a pair of integers such that $K^2=2\chi-5$ and $\chi\geq 3$.
Let us assume that $K^2+1\in 8\cdot\mathbb{Z}$ or $(\chi,K^2)=(7,9)$.
The connected components of $\mathfrak{M}_{\chi,K^2}$, that we denote by 
$\mathfrak{M}_{\chi,K^2}^{I}$ and $\mathfrak{M}_{\chi,K^2}^{II}$, will be described below.
%We will denote by $\mathfrak{M}_{\chi,K^2}^{I}$ (resp. $\mathfrak{M}_{\chi,K^2}^{II}$) 
%the connected component of $\mathfrak{M}_{\chi,K^2}$ containing  the surfaces on the first
%(resp. second) deformation class.
\end{remark}

In the case $K^2=2\chi-6$ the canonical image of a surface completely determines its deformation class. 
The following results show that in the case $K^2=2\chi-5$ sometimes the canonical image of a surface 
$X$ is not enough to determine its deformation class. When $(\chi,K^2)=(6,7)$ or $(\chi,K^2)=(7,9)$
we may need to know $X'$ too.

\begin{theorem}[{{\cite[Theorem 8.1]{Hor2}}}]
\label{DefClass2Gen}
Let $(\chi,K^2)$ be a pair of integers such that $K^2=2\chi-5$, $\chi\geq 3$, $K^2+1\in 8\cdot\mathbb{Z}$ 
and $K^2\neq 7$. The image of the canonical map of a smooth minimal surface in 
$\mathfrak{M}_{\chi,K^2}^{I}$ is $\mathbb{F}_e$ for some $e\in\{1,3,\ldots,\frac{1}{4}(K^2+1)-1\}$.
The image of the canonical map of a smooth minimal surface in $\mathfrak{M}_{\chi,K^2}^{II}$ is 
$\mathbb{F}_{\frac{1}{4}(K^2+1)+1}$.
\end{theorem}

\begin{theorem}[{{\cite[Theorem 6.1]{Hor2}}}]
\label{DefClass267}
Let $(\chi,K^2)=(6,7)$. Smooth minimal surfaces $X$ in $\mathfrak{M}_{\chi,K^2}^{I}$ have canonical
image $\mathbb{F}_1$ and $X'$ is a double cover of a blow-up of $\mathbb{F}_1$ at two points of a fiber. 
Smooth minimal surfaces $X$ in $\mathfrak{M}_{\chi,K^2}^{II}$ either have canonical image a cone over a rational 
curve of degree $3$ in $\mathbb{P}^3$ or they have canonical image $\mathbb{F}_1$ and $X'$ is a double cover of 
a blow-up of $\mathbb{F}_1$ at a point on the negative section.
\end{theorem}

\begin{theorem}[{{\cite[Theorem 7.1]{Hor2}}}]
\label{DefClass279}
Let $(\chi,K^2)=(7,9)$. Smooth minimal surfaces $X$ in $\mathfrak{M}_{\chi,K^2}^{I}$ have canonical image
$\mathbb{F}_e$ for some $e\in\{0,2\}$ and $X'$ is a double cover of a blow-up of $\mathbb{F}_e$ at 
two points of a fiber. Smooth minimal surfaces $X$ in $\mathfrak{M}_{\chi,K^2}^{II}$ have canonical image
$\mathbb{F}_2$ and $X'$ is a double cover of $\mathbb{F}_2$.
\end{theorem}

\section{$\mathbb{Z}_2^2$-covers.}\label{Z22 covers.}
A $\mathbb{Z}_2^2$-cover of a variety $Y$ is a finite map $f:X\rightarrow Y$ together with a faithful action of 
$\mathbb{Z}_2^2$ on $X$ such that $f$ exhibits $Y$ as $X/\mathbb{Z}_2^2$.
The structure theorem for smooth $\mathbb{Z}_2^2$-covers was first 
given by Catanese \cite{Cata1}.  
According to {\cite[Section 2]{Cata1}} or {\cite[Theorem 2.1]{Par1}},
to define a smooth $\mathbb{Z}_2^2$-cover of a smooth variety $Y$ it suffices to consider both:
\begin{enumerate}
 \item[-] Smooth divisors $D_1, D_2, D_3$ such that the branch locus $B=D_1+D_2+D_3$ is a normal crossing divisor.
 \item[-] Line bundles $L_1, L_2, L_3$ satisfying $2L_1\equiv D_2+D_3, 2L_2\equiv D_1+D_3$ and such that  $L_3\equiv L_1+L_2-D_3$.
\end{enumerate}
The set $\{L_i,D_j\}_{i,j}$ is called the building data of the cover.

\begin{remark}\label{Z22sing}
 We can also consider $\mathbb{Z}_2^2$-covers of a variety $Y$ with singularities. If $Y$ is normal 
 similar building data is required to define the cover, but there are some differences 
 (see \cite{Alexeevpardini1}):
 \begin{enumerate}
 \item[-] The divisors $D_j$ are not necessarily Cartier divisors but Weil divisors.
 \item[-] The sheaves $L_i$ are not necessarily invertible sheaves but reflexive divisorial sheaves.
\end{enumerate}
In addition, the branch locus may have non-divisorial components. 
\end{remark}

\begin{proposition}
[{{\cite[Section 2]{Cata1}}} or {{\cite[Proposition 4.2]{Par1}}}] 
Let $Y$ be a smooth surface and $f:X\rightarrow Y$ a smooth $\mathbb{Z}_2^2$-cover with building data 
$\{L_i,D_j\}_{i,j}$. Then:
\begin{equation*}
\begin{split}
2K_X\equiv f^*(2K_Y+D_1+D_2+D_3),\\
K_X^2=(2K_Y+D_1+D_2+D_3)^2,\\
p_g(X)=p_g(Y)+\sum_{i=1}^3h^0(K_Y+L_i),\\
\chi(\mathcal{O}_X)=4\chi(\mathcal{O}_Y)+\frac{1}{2}\sum_{i=1}^3L_i(L_i+K_Y).
\end{split}
\end{equation*}
\end{proposition}

\begin{remark}\label{canonicalimage}
 Let us consider one of the line bundles defining a $\mathbb{Z}^2_2$-cover $f:X\rightarrow Y$, for instance $L_1$.
 Then we can decompose $f$ as $f_2\circ f_1$ where $f_2:X_1\rightarrow Y$ is a double cover of $Y$ branched along 
 $D_2+D_3$ and $f_1:X\rightarrow X_1$ is a double cover of $X_1$ branched along 
 $f_2^*D_{1}\cup f_2^*(D_{2}\cap D_{3})$. %If $D_{2}\cap D_{3}=\varnothing$ then 
 In particular $K_X=f_1^*(K_{X_1}+f_2^*L_1)$. Let us denote $N= h^0(K_{X_1}+f_2^*L_1)$ and 
 $i:X_1\rightarrow \mathbb{P}^{N-1}$ the morphism defined by the complete linear system $|K_{X_1}+f_2^*L_1|$.
 Then $(i\circ f_1)^*(\mathcal{O}_{\mathbb{P}^{N-1}}(1))=K_X$. Hence $i\circ f_1$ is a 
 morphism induced by some sections of $K_X$. It follows that if $h^0(K_X)=h^0( K_{X_1}+f_2^*L_1)$ 
 then $i\circ f_1$ is the canonical map of $X$ and $i(X_1)$ is its canonical image.   
\end{remark}

\section{Constructions when $K^2=2\chi-6$.}\label{Constructions.}
In this section we are going to prove Theorem \ref{MainResult} in the case $K^2=2\chi-6$ as follows.
We fix a pair of integers $(\chi,K^2)$ such that $K^2=2\chi-6$ and $\chi\geq 4$.
If $\chi$ is even then $K^2\notin 8\cdot \mathbb{Z}$ and $\mathfrak{M}_{\chi,K^2}$ has a unique connected 
component by Theorem \ref{DefClass}. Then we just need to find a surface in $\mathfrak{M}_{\chi,K^2}$ with
a $\mathbb{Z}^2_2$-action.
If $\chi$ is odd it may happen that $K^2\in 8\cdot \mathbb{Z}$ and in this case $\mathfrak{M}_{\chi,K^2}$ 
has two connected components again by Theorem \ref{DefClass}. What we are going to do is to find a 
surface in $\mathfrak{M}_{\chi,K^2}$ and by studying its canonical image we are going to check that it
does not belong to the second deformation class when $K^2\in 8\cdot \mathbb{Z}$. Finally, we are going 
to find surfaces in the second deformation class
for any $K^2\in 8\cdot \mathbb{Z}$.

Let us assume first that $\chi$ is even. 
We denote by $\mathbb{F}_2$ the Hirzebruch surface with negative section $\Delta_0$ of self-intersection $(-2)$
and fiber $F$. The smooth $\mathbb{Z}_2^2$-cover $\pi:S\rightarrow \mathbb{F}_2$ of $\mathbb{F}_2$ whose branch 
locus $B=D_1+D_2+D_3$ is a normal crossing divisor consisting of smooth and irreducible divisors
\begin{equation*}
\begin{split}
D_1\in |\Delta_0|\\
D_2\in |\Delta_0+2 F|\\
D_3\in |3\Delta_0+ (\chi + 4)F|
\end{split}
\end{equation*}
satisfies $\chi(\mathcal{O}_S)=\chi$ and $K^2_S=K^2$. The canonical divisor of $S$ is big and nef because 
$2K_S$ is the pullback via $\pi$ of the big and nef divisor $\Delta_0+(\chi-2)F$. Hence $S$ is minimal and its canonical model, that has a 
$\mathbb{Z}^2_2$-action, belongs to $\mathfrak{M}_{\chi,K^2}$. We notice that when $\chi>4$ the divisor $K_S$ is ample and $S$ itself is  a canonical model. When $\chi=4$ it follows from Remark \ref{Case42} that $K_S$ is not ample because $S$ has a genus $2$ fibration.
The canonical model of $S$, obtained by contracting the $(-2)$-curve $(\pi^*\Delta_0)_{\text{red}}$,
is a $\mathbb{Z}_2^2$-cover of the quadric cone.

Let us assume now that $\chi$ is odd, which happens always if $K^2\in 8\cdot \mathbb{Z}$. 
We denote by $\Delta_0$ and $F$ the two classes of fibers of $\mathbb{P}^1\times \mathbb{P}^1$. 
The smooth $\mathbb{Z}_2^2$-cover $\pi:S\rightarrow \mathbb{P}^1\times \mathbb{P}^1$ of 
$\mathbb{P}^1\times \mathbb{P}^1$ whose branch locus $B=D_1+D_2+D_3$ is a normal crossing divisor
consisting of smooth and irreducible divisors
\begin{equation*}
\begin{split}
D_1\in |\Delta_0|\\
D_2\in |\Delta_0|\\
D_3\in |3\Delta_0+ (\chi + 1)F|
\end{split}
\end{equation*}
satisfies $\chi(\mathcal{O}_S)=\chi$ and $K^2_S=K^2$. The canonical divisor of $S$ is ample because 
$2K_S$ is the pullback via $\pi$ of the ample divisor $\Delta_0+(\chi-3)F$. Hence
$S\in \mathfrak{M}_{\chi,K^2}$ and it has a 
$\mathbb{Z}^2_2$-action. According to Remark \ref{canonicalimage} we can decompose
$\pi$ as $\pi=\pi_2\circ\pi_1$ where 
$\pi_2:\mathbb{P}^1\times \mathbb{P}^1\rightarrow \mathbb{P}^1\times \mathbb{P}^1$ is a double cover of 
$\mathbb{P}^1\times \mathbb{P}^1$ branched along $D_1+D_2$ and 
$\pi_1:S\rightarrow \mathbb{P}^1\times \mathbb{P}^1$
is a double cover of $\mathbb{P}^1\times \mathbb{P}^1$ branched along 
$\pi_2^*D_3\in|6\Delta_0+(\chi +1)F|$. Then $K_S$ is the pullback via $\pi_2$ of the very ample divisor 
$K_{\mathbb{P}^1\times \mathbb{P}^1}+\frac{1}{2}\pi_2^*D_3\in|\Delta_0+(\frac{\chi-3}{2})F|$. Since 
$h^0(K_S)=h^0(\mathbb{P}^1\times \mathbb{P}^1, \Delta_0+(\frac{\chi-3}{2})F)$
the canonical image of $S$ is $\mathbb{P}^1\times \mathbb{P}^1$. When $K^2\in 8\cdot \mathbb{Z}$
the surface $S$ belongs to the first deformation class  by Theorem \ref{DefClass}.

Finally, we are going to consider a smooth $\mathbb{Z}_2^2$-cover $\pi:S\rightarrow \mathbb{F}_e$ of the Hirzebruch surface
$\mathbb{F}_e$ with negative section $\Delta_0$ of self-intersection $(-e)\leq -2$ and fiber $F$. 
In this case, the branch locus $B=D_1+D_2+D_3$ is a normal crossing divisor consisting of
\begin{align*}
& D_1\in |F|& & D_1=0\nonumber\\
& D_2\in |F| \qquad\qquad\qquad \text{if $e\notin 2\cdot\mathbb{Z}$}& &D_2\in |2F|\qquad\qquad\qquad 
\text{if $e\in 2\cdot\mathbb{Z}$}\\
& D_3\in |6\Delta_0+5eF| & &D_3\in |6\Delta_0+5eF|
\end{align*}
where $D_3$ is the union of the negative section $\Delta_0$ and a smooth and irreducible divisor 
$D_3'\in |5\Delta_0+5eF|$. Then $\chi(\mathcal{O}_S)=4e-1$ and $K_S^2=8e-8$. 
In particular $K^2_S=2\chi(\mathcal{O}_S)-6$ and $K_S^2\in 8\cdot\mathbb{Z}$.
According to Remark \ref{canonicalimage} we can decompose $\pi$ as $\pi=\pi_2\circ\pi_1$ 
where $\pi_2:\mathbb{F}_{2e}\rightarrow \mathbb{F}_e$ is the double cover of $\mathbb{F}_e$ 
branched along $D_1+D_2$ and $\pi_1:S\rightarrow \mathbb{F}_{2e}$ is the double cover of the Hirzebruch
surface $\mathbb{F}_{2e}$ with negative section $\Gamma_0$ of self-intersection $(-2e)$ and fiber 
$G$ branched along 
$\pi_2^*D_3\in|6\Gamma_0+10eG|$. Then $K_S$ is the pullback via $\pi_2$ of the divisor 
$K_{\mathbb{F}_{2e}}+\frac{1}{2}\pi_2^*D_3\in|\Gamma_0+(3e-2)G|$. Since 
$h^0(K_S)=h^0(\mathbb{F}_{2e}, \Gamma_0+(3e-2)G)$ the canonical image of $S$ is the image of 
$\mathbb{F}_{2e}$ via the map $\varphi$ defined by the complete linear system $|\Gamma_0+(3e-2)G|$. 
If $e>2$ this divisor is very ample and therefore $\varphi$ is an embedding. 
If $e=2$ then $\varphi$ embeds $\mathbb{F}_{2e}$ as a surface of minimal
degree.
In any case the surface $S$ belongs to the second deformation class by Theorem 
\ref{DefClass}.

\section{Constructions when $K^2=2\chi-5$.}\label{ConstructionsSecondCase.}
In this section we are going to prove Theorem \ref{MainResult} in the case $K^2=2\chi-5$ as follows.
We fix a pair of integers $(\chi,K^2)\neq (7,9)$ such that $K^2=2\chi-5$ and $\chi\geq 3$.
If $\chi$ is odd then $K^2+1\notin 8\cdot \mathbb{Z}$ and $\mathfrak{M}_{\chi,K^2}$ has a unique connected 
component by Theorem \ref{DefClass2}. Then we just need to find a surface in $\mathfrak{M}_{\chi,K^2}$ with 
a $\mathbb{Z}^2_2$-action.
If $\chi$ is even it may happen that $K^2+1\in 8\cdot \mathbb{Z}$ and in this case $\mathfrak{M}_{\chi,K^2}$ 
has two connected components again by Theorem \ref{DefClass2}. What we are going to do is to find a surface
in $\mathfrak{M}_{\chi,K^2}$ and by studying its canonical image we are going to check that it does not belong 
to $\mathfrak{M}_{\chi,K^2}^{II}$ when $K^2+1\in 8\cdot \mathbb{Z}$. Finally, we are going to find surfaces in 
$\mathfrak{M}_{\chi,K^2}^{II}$ when $K^2+1\in 8\cdot \mathbb{Z}$.
The case $(\chi,K^2)=(7,9)$ will be studied separately.

Let us assume first that $\chi$ is odd and bigger than $3$. We denote by $\mathbb{F}_3$ the Hirzebruch surface
with negative section $\Delta_0$ of self-intersection $(-3)$
and fiber $F$. Then the smooth $\mathbb{Z}_2^2$-cover $\pi:S\rightarrow \mathbb{F}_3$ of $\mathbb{F}_3$ 
whose branch 
locus $B=D_1+D_2+D_3$ is a normal crossing divisor consisting of smooth and irreducible divisors
\begin{equation*}
\begin{split}
D_1\in |\Delta_0|\\
D_2\in |\Delta_0+4 F|\\
D_3\in |3\Delta_0+ (\chi + 5)F|
\end{split}
\end{equation*}
satisfies $\chi(\mathcal{O}_S)=\chi$ and $K^2_S=K^2$. The canonical divisor of $S$ is ample because 
$2K_S$ is the pullback via $\pi$ of the ample divisor $\Delta_0+(\chi-1)F$. Therefore,
 $S\in \mathfrak{M}_{\chi,K^2}$ and it has a $\mathbb{Z}^2_2$-action. 
 Let us suppose that $\chi=3$. Then we cannot assume $D_3$ to be smooth in the previous example. Nevertheless
 a smooth $\mathbb{Z}_2^2$-cover $\pi:S\rightarrow \mathbb{P}^2$ of $\mathbb{P}^2$ whose branch 
locus $B=D_1+D_2+D_3$ is a normal crossing divisor consisting of smooth and irreducible divisors
\begin{equation*}
\begin{split}
D_1\in |\mathcal{O}_{\mathbb{P}^2}(1)|\\
D_2\in |\mathcal{O}_{\mathbb{P}^2}(1)|\\
D_3\in |\mathcal{O}_{\mathbb{P}^2}(5)|
\end{split}
\end{equation*}
satisfies $(\chi(\mathcal{O}_S), K^2_S)=(3,1)$. The canonical divisor of $S$ is ample because 
$2K_S$ is the pullback via $\pi$ of the ample sheaf $\mathcal{O}_{\mathbb{P}^2}(1)$. Hence
 $S\in \mathfrak{M}_{3,1}$ and it has a $\mathbb{Z}^2_2$-action.
 
Let us assume now that $\chi$ is even and greater or equal than $4$, which happens always if 
$K^2+1\in 8\cdot \mathbb{Z}$. We denote by $\mathbb{F}_1$ the Hirzebruch surface with negative
section $\Delta_0$ of self-intersection $(-1)$ and fiber $F$. 
Then the smooth $\mathbb{Z}_2^2$-cover $\pi:S\rightarrow \mathbb{F}_1$ of $\mathbb{F}_1$ whose branch 
locus $B=D_1+D_2+D_3$ is a normal crossing divisor consisting of smooth and irreducible divisors
\begin{equation*}
\begin{split}
D_1\in |\Delta_0|\\
D_2\in |\Delta_0+2 F|\\
D_3\in |3\Delta_0+ (\chi + 2)F|
\end{split}
\end{equation*}
satisfies $\chi(\mathcal{O}_S)=\chi$ and $K^2_S=K^2$. The canonical divisor of $S$ is ample because 
$2K_S$ is the pullback via $\pi$ of the ample divisor $\Delta_0+(\chi-2)F$. 
Hence $S\in \mathfrak{M}_{\chi,K^2}$ and it has a 
$\mathbb{Z}^2_2$-action. Let us restrict now to the case $K^2+1\in 8\cdot \mathbb{Z}$. 
If we denote by $q:\widetilde{\mathbb{F}}_2\rightarrow \mathbb{F}_1$
 the blow-up  
 of $\mathbb{F}_1$ at the point of intersection of $D_1$ and $D_2$ with exceptional divisor $E'$,
the cover $\pi$
induces a
 $\mathbb{Z}_2^2$-cover $\widetilde{\pi}:\widetilde{S}\rightarrow \widetilde{\mathbb{F}}_2$ with
 branch locus $\widetilde{B}=\widetilde{D}_1+\widetilde{D}_2+\widetilde{D}_3$ where 
 \begin{equation*}
\begin{split}
\widetilde{D}_1=q^*D_1-E'\\
\widetilde{D}_2=q^*D_2-E'\\
\widetilde{D}_3=q^*D_3+E'
\end{split}
\end{equation*} 
If we denote by $\pi_2:\widetilde{\mathbb{F}}_1\rightarrow \widetilde{\mathbb{F}}_2$
the intermediate $\mathbb{Z}_2$-cover of $\widetilde{\pi}$ 
with branch locus $\widetilde{D}_1+\widetilde{D}_2$ (see Remark \ref{canonicalimage}), we have that
$\widetilde{\mathbb{F}}_1$
is a blow-up 
 $b:\widetilde{\mathbb{F}}_1\rightarrow \mathbb{F}_1$
 of $\mathbb{F}_1$ at two points of a fiber. In addition the induced map 
 $\varepsilon:\widetilde{S}\rightarrow S$ is the blow-up of $S$ at the point over $D_1\cap D_2$
 with exceptional divisor $E$. We obtain the following commutative diagram:
 \begin{equation*}
 \xymatrix@R-1pc{
S
%\ar@{-->}@/^2pc/[rr]^{\varphi}
\ar[dd]_{\pi}&\widetilde{S}\ar[l]_{\varepsilon}\ar[d]_{\pi_1}\ar[r]^{f} &\mathbb{F}_1\\
&\widetilde{\mathbb{F}}_1\ar[d]_{\pi_2}\ar[ur]_{b} & \\
\mathbb{F}_1&\widetilde{\mathbb{F}}_2\ar[l]^{q} & 
}
\end{equation*}
where $\pi_1$ is the $\mathbb{Z}_2$-cover of $\widetilde{\mathbb{F}}_1$ branched along $\pi_2^*\widetilde{D}_3$ and 
%$\varphi$ is the rational map induced by 
$f:=b\circ \pi_1$.
 I claim that we are in case $a)$ of Theorem 
\ref{CanSys5} and the canonical image of $S$ is $\mathbb{F}_1$. Indeed, let us denote by 
$\Sigma_0$ (resp. $G$) the negative section (resp. the class of a fiber) of the latter Hirzebruch surface 
$\mathbb{F}_1$. Taking into account the standard formulas for $\mathbb{Z}_2$-covers and the fact that 
$K_{\widetilde{S}}=\varepsilon^*K_S+E$, a straightforward calculation gives:
\begin{equation}\label{CanImageSupport}
 \varepsilon^*K_S=f^*\left(\Sigma_0+\frac{\chi-2}{2}G\right)+E.
\end{equation}
Now, 
it follows from the equality $h^0(K_S)=h^0(\Sigma_0+\frac{\chi-2}{2}G)$ together with (\ref{CanImageSupport}) that $E$ is the fixed part of the linear system $|\varepsilon^*K_S|$ and $|f^*\left(\Sigma_0+\frac{\chi-2}{2}G\right)|$ is its moving part. We conclude that 
the canonical image of $S$ coincides with the image of the map induced by the complete linear system 
$|f^*(\Sigma_0+\frac{\chi-2}{2}G)|$. Since $\Sigma_0+\frac{\chi-2}{2}G$ is very ample (recall the assumption $K^2+1\in 8\cdot \mathbb{Z}$ that implies $\chi\geq 6$) we conclude that the canonical image of $S$ is $\mathbb{F}_1$ as claimed.
Therefore $S\in \mathfrak{M}_{\chi,K^2}^{I}$ by Theorem \ref{DefClass2Gen}.

Now we are going to consider a $\mathbb{Z}_2^2$-cover $\pi:T\rightarrow \mathbb{F}_{k+1}$ of the Hirzebruch surface
$\mathbb{F}_{k+1}$ with negative section $\Delta_0$ of self-intersection $-(k+1)\leq -3$ and fiber $F$. 
In this case, the branch locus $B=D_1+D_2+D_3$ consists of
\begin{align*}
& D_1\in |F|& & D_1=0\nonumber\\
& D_2\in |F| \qquad\qquad\qquad\qquad \text{if $k\in 2\cdot\mathbb{Z}$}& &D_2\in |2F|\qquad\qquad\qquad\qquad 
\text{if $k\notin 2\cdot\mathbb{Z}$}\\
& D_3\in |6\Delta_0+5(k+1)F| & &D_3\in |6\Delta_0+5(k+1)F|
\end{align*}
The divisor $D_3$ is the union of the negative section $\Delta_0$ and an irreducible divisor 
$D_3'\in |5\Delta_0+5(k+1)F|$ with an ordinary triple point $p$. The divisor $D_2$ passes through 
$p$ intersecting $D_3$ with multiplicity $3$. 
Then $\chi(\mathcal{O}_T)=4k+3$ and $K_T^2=8k$. Nevertheless $T$ has an elliptic singularity over $p$ such that
its minimal resolution $r:S\rightarrow T$ satisfies $\chi(\mathcal{O}_S)=4k+2$ and
$K_S^2=8k-1$ (see Remark \ref{EllipticSingRemark} below). In particular $K_S^2=2\chi(\mathcal{O}_S)-5$ and $K_S^2+1\in 8\cdot\mathbb{Z}$.
If we denote 
$q_2:\widetilde{\mathbb{F}}_{k+1}\rightarrow \mathbb{F}_{k+1}$ the blow-up  of $\mathbb{F}_{k+1}$ at $p$ with
exceptional divisor $E_1'$,
the cover $\pi$ induces another $\mathbb{Z}_2^2$-cover $\widetilde{\pi}:S\rightarrow \widetilde{\mathbb{F}}_{k+1}$.
The canonical divisor of $S$ is ample because $2K_S$ is the pull-back via $\widetilde{\pi}$
of $q_2^*(2\Delta_0+(3k+1)F)-E'_1=:D$ and this divisor is ample by the Nakai-Moishezon criterion.
Indeed, let us suppose that there exists an irreducible curve $C\in |q_2^*(a\Delta_0+bF)-cE'_1|$ such that $CD<0$ for some non negative integers $a,b,c$. This implies that $q_2(C)\in |a\Delta_0+bF|$ is an irreducible 
curve of $\mathbb{F}_{k+1}$ with a point of multiplicity $c>(k-1)a+2b$, which is clearly impossible.
Now, 
$\pi$ induces one more $\mathbb{Z}_2^2$-cover $\widetilde{\widetilde{\pi}}:\widetilde{S}\rightarrow 
\widetilde{\widetilde{\mathbb{F}}}_{k+1}$ on the blow-up 
$\widetilde{\widetilde{\mathbb{F}}}_{k+1}\rightarrow \widetilde{\mathbb{F}}_{k+1}$ of 
$\widetilde{\mathbb{F}}_{k+1}$ at 
the point of intersection $p'$ of $E_1'$ with the strict transform via $q_2$ of the fiber of 
$\mathbb{F}_{k+1}$ through $p$.
It turns out that one of the 
 intermediate surfaces of the cover $\widetilde{\widetilde{\pi}}$ is a blow-up 
 $\widetilde{\mathbb{F}}_{2k+2}\rightarrow \mathbb{F}_{2k+1}$
 of $\mathbb{F}_{2k+1}$ at two points of a fiber, one of them belonging to the negative section and
 $\widetilde{S}\rightarrow S$ is the blow-up of $S$ at the point over $p'$.
It can be proved as above that we are in 
case $a)$ of Theorem \ref{CanSys5} and the canonical image of $S$ is $\mathbb{F}_{2k+1}$.
Therefore $S\in \mathfrak{M}_{\chi,K^2}^{II}$ by Theorem \ref{DefClass2Gen}
and it has a $\mathbb{Z}_2^2$-action.

When $k=1$, i.e. when $(\chi, K^2)=(6,7)$, we can consider the same covers.
Nevertheless, the canonical divisor of $S$ is not ample because $(\widetilde{\pi}^*q_2^*\Delta_0)_{\text{red}}$ 
is a $(-2)$-curve.
Its canonical model $Z$, obtained by contracting this $(-2)$-curve, is a $\mathbb{Z}_2^2$-cover of a blow-up
of the quadric cone.
Arguing as before, and taking into account Remark \ref{Z22sing}, it can be proved that the canonical divisor
of $Z$ is ample and its canonical image is 
a cone over a rational curve of degree $3$ in $\mathbb{P}^3$. Hence $Z\in\mathfrak{M}_{6,7}^{II}$ by 
Theorem \ref{DefClass267} and it has a $\mathbb{Z}_2^2$-action.

Let us study now the case $(\chi,K^2)=(7,9)$. To obtain a smooth minimal surface in $\mathfrak{M}_{\chi,K^2}^{I}$
it suffices to consider the example that we constructed in the general case with $\chi$ odd. Indeed, 
the smooth $\mathbb{Z}_2^2$-cover $\pi:S\rightarrow \mathbb{F}_3$ of $\mathbb{F}_3$ whose branch 
locus $B=D_1+D_2+D_3$ is a normal crossing divisor consisting of smooth and irreducible divisors
\begin{equation*}
\begin{split}
D_1\in |\Delta_0|\\
D_2\in |\Delta_0+4 F|\\
D_3\in |3\Delta_0+12F|
\end{split}
\end{equation*}
satisfies $(\chi(\mathcal{O}_S), K^2_S)=(7,9)$. The canonical divisor of $S$ is ample because 
$2K_S$ is the pullback via $\pi$ of the ample divisor $\Delta_0+6F$. The cover $\pi$ induces another
 $\mathbb{Z}_2^2$-cover $\widetilde{S}\rightarrow \widetilde{\mathbb{F}}_4$
 where 
 $\widetilde{\mathbb{F}}_4\rightarrow \mathbb{F}_3$ 
 is the blow-up  
 of $\mathbb{F}_3$ at the point of intersection of $D_1$ and $D_2$. It turns out that
 one of the 
 intermediate surfaces of this $\mathbb{Z}_2^2$-cover is a blow-up 
 $\widetilde{\mathbb{F}}_2\rightarrow \mathbb{F}_2$
 of $\mathbb{F}_2$ at two points of a fiber and
 $\widetilde{S}\rightarrow S$ is the blow-up of $S$ at the point over $D_1\cap D_2$. 
 It can be proved as above that we are in case $a)$ of Theorem 
\ref{CanSys5} and the canonical image of $S$ is $\mathbb{F}_2$.
Therefore $S\in \mathfrak{M}_{\chi,K^2}^{I}$ by Theorem \ref{DefClass279}
and it has a $\mathbb{Z}_2^2$-action.

To obtain a smooth minimal surface in $\mathfrak{M}_{\chi,K^2}^{II}$ we can consider 
a smooth $\mathbb{Z}_2^2$-cover $\pi:S\rightarrow \mathbb{P}^2$ of $\mathbb{P}^2$ whose branch 
locus $B=D_1+D_2+D_3$ is a normal crossing divisor consisting of smooth and irreducible divisors
\begin{equation*}
\begin{split}
D_1\in |\mathcal{O}_{\mathbb{P}^2}(1)|\\
D_2\in |\mathcal{O}_{\mathbb{P}^2}(1)|\\
D_3\in |\mathcal{O}_{\mathbb{P}^2}(7)|
\end{split}
\end{equation*}
This cover satisfies $(\chi(\mathcal{O}_S), K^2_S)=(7,9)$ and the canonical divisor of $S$ is ample because 
$2K_S$ is the pullback via $\pi$ of the ample sheaf $\mathcal{O}_{\mathbb{P}^2}(3)$. 
In addition, $\pi$ induces another $\mathbb{Z}_2^2$-cover $\widetilde{S}\rightarrow \mathbb{F}_1$
where $\mathbb{F}_1\rightarrow \mathbb{P}^2$ is the blow-up
of $\mathbb{P}^2$ at the point $p$ of intersection of $D_1$
and $D_2$.
It turns out that one of the intermediate surfaces of this $\mathbb{Z}_2^2$-cover is $\mathbb{F}_2$
and
 $\widetilde{S}\rightarrow S$ is the blow-up of $S$ at the point over $D_1\cap D_2$.
Arguing as usual it can be proved that 
we are in case 
$c)$ of Theorem \ref{CanSys5} and
the canonical image of $S$ is $\mathbb{F}_2$. Therefore 
$S\in \mathfrak{M}_{\chi,K^2}^{II}$ by Theorem \ref{DefClass279}
and it has a $\mathbb{Z}_2^2$-action.

\begin{remark}\label{EllipticSingRemark}
 Let us fix an admissible pair $(\chi, K^2)$ such that 
$K^2=2\chi-5$, $\chi\geq 7$, $K^2+1\in 8\cdot \mathbb{Z}$. Our examples of surfaces $X\in 
\mathfrak{M}_{\chi, K^2}^{II}$ are constructed as minimal resolutions of surfaces $Y$ with a 
 singularity obtained by contracting a smooth elliptic curve with self-intersection $(-1)$.

Our procedure to obtain these singularities (see \cite{Cata2}) 
consists in constructing a $\mathbb{Z}_{2}^2$-cover $Y\rightarrow \mathbb{F}_e$ with branch locus
$B=D_1+D_2+D_3$ such that $D_3$ has an ordinary triple point $p$, $D_2$ passes through $p$ intersecting $D_3$ with multiplicity $3$ and $D_1$ does not pass through $p$.  We can resolve the singularity over $p$ in a canonical way. Let $b:\widetilde{\mathbb{F}}_e\rightarrow \mathbb{F}_e$ be the blow-up of $\mathbb{F}_e$ at $p$ with exceptional divisor $E$. Then there is 
 a $\mathbb{Z}_2^2$-cover $X\rightarrow\widetilde{\mathbb{F}}_e$ with branch locus $\widetilde{B}=\widetilde{D_1}+\widetilde{D_2}+
\widetilde{D_3}$ where 
\begin{equation*}
 \begin{split}
  \widetilde{D_1}=b^*D_1-3E\\
  \widetilde{D_2}=b^*D_2-E\\
  \widetilde{D_3}=b^*D_3+E
 \end{split}
\end{equation*}
The induced map $X\rightarrow Y$ resolves the singularity over $p$.
In addition, using the formulas for $\mathbb{Z}_2^2$-covers we can check that the
holomorphic Euler characteristic and the self intersection of the canonical class of $X$ are
reduced in $1$ with respect to those of $Y$.  In particular $K^2_Y=2\chi(\mathcal{O}_Y)-6$. 
Nevertheless, the surfaces $Y$ so constructed do not belong to 
the moduli space $\mathfrak{M}_{\chi+1, K^2+1}$
because these elliptic singularities are not canonical.
They belong to the KSBA-compactification $\overline{\mathfrak{M}}_{\chi+1, K^2+1}$ of
$\mathfrak{M}_{\chi+1, K^2+1}$.

A question which arises now is whether they are in the closure of
$\mathfrak{M}_{\chi+1, K^2+1}$ or if, on the contrary, they belong to 
 a different irreducible component of $\overline{\mathfrak{M}}_{\chi+1, K^2+1}$.
  To answer this question we notice that 
  we can assume $D_3$ to be smooth
 and the surface $Z$ so obtained has the same invariants as $Y$ but is not 
 singular.
 Since $Y$ and $Z$ are $\mathbb{Z}_2^2$-covers of the same surface  and the divisors
 forming the branch locus belong to the same linear systems in both cases, it follows that 
 $Y\in \overline{\mathfrak{M}}_{\chi+1, K^2+1}$ is a natural deformation (see \cite[Section 5]{Par1})
of $Z\in \mathfrak{M}_{\chi+1, K^2+1}$. This shows that $Y$ belongs to the closure of
$\mathfrak{M}_{\chi+1, K^2+1}$ in $\overline{\mathfrak{M}}_{\chi+1, K^2+1}$.
\end{remark}

\section{Genus $2$ fibrations.} \label{Genus2Fibrations}
In \cite{HorGen2} Horikawa studied genus $2$ fibrations. One of the things he showed is that if we have a genus $2$
pencil $f:X\rightarrow \mathbb{P}^1$, then the following formula holds:
\begin{equation*}
 K^2_X=2\chi(\mathcal{O}_X)-6+\sum_{p\in\mathbb{P}^1}H(p).
\end{equation*}
The contribution $H(p)$ is a non-negative integer that is bigger than $0$ if and only if 
the fiber over $p\in\mathbb{P}^1$ is not $2$-connected.

Most of our examples of surfaces on the line 
$K^2=2\chi-6$ are genus $2$ pencils (see Proposition \ref{Gen26}) such that all their fibers are $2$-connected.
Most of our examples of surfaces on the line 
$K^2=2\chi-5$ are genus $2$ pencils (see Proposition \ref{Gen25}) with a unique not $2$-connected 
fiber with contribution $H(p)=1$.
Moreover, this fiber is usually obtained by making the divisors
$D_1$ and $D_2$ of the branch locus of the $\mathbb{Z}_2^2$-cover intersect and consists
of two elliptic curves with self-intersection 
$(-1)$ intersecting transversally in one point.

\noindent \begin{acknowledgements}
The author is deeply indebted to his supervisor Margarida Mendes Lopes for all her help.
The author also thanks the anonymous reviewer for her/his thorough reading of the paper and suggestions.
\end{acknowledgements}

\bibliographystyle{plain}      
\bibliography{Z22ActionsOnHorikawaSurfaces}

\begin{thebibliography}{10}

\bibitem{Alexeevpardini1}
V.~Alexeev and R.~Pardini.
\newblock Non-normal abelian covers.
\newblock {\em Compositio Mathematica}, 148(4):1051--1084, 2012.

\bibitem{Barth}
W.P. Barth, K.~Hulek, C.A.M. Peters, and A.~Van~de Ven.
\newblock Compact complex surfaces.
\newblock {\em Ergebnisse der Mathematik und ihrer Grenzgebiete. 3. Folge. A
  $\text{Series}$ of Modern Surveys in Mathematics}, 4, 2004.

\bibitem{Bomb}
E.~Bombieri.
\newblock Canonical models of surfaces of general $\text{type}$.
\newblock {\em Publ. Math. IHES}, 42:171–219, 1973.

\bibitem{Cata1}
F.~Catanese.
\newblock On the moduli spaces of surfaces of general $\text{Type}$.
\newblock {\em J. Diff. Geom.}, 19(2):483--515, 1984.

\bibitem{Cata2}
F.~Catanese.
\newblock Singular bidouble covers and the construction of interesting
  algebraic surfaces.
\newblock {\em Contemporary Mathematics}, 241:97--120, 1999.

\bibitem{CataDef}
F.~Catanese.
\newblock Differentiable and deformation $\text{type}$ of algebraic surfaces,
  real and symplectic structures.
\newblock {\em Springer lecture notes in mathematics 1938}, Symplectic
  4-manifolds and algebraic surfaces:55–167, 2008.

\bibitem{Hor1}
E.~Horikawa.
\newblock Algebraic surfaces of general $\text{type}$ with small $c_1^2$,
  $\text{I}$.
\newblock {\em Ann. Math.}, 104:357--387, 1976.

\bibitem{Hor2}
E.~Horikawa.
\newblock Algebraic surfaces of general $\text{type}$ with small $c_1^2$,
  $\text{II}$.
\newblock {\em Inventiones matematicae}, 37:121–155, 1976.

\bibitem{HorGen2}
E.~Horikawa.
\newblock On algebraic surfaces with pencils of curves of genus $2$.
\newblock {\em Cambridge University Press}, pages 79--90, 1977.

\bibitem{Mendespardini1}
M.M. Lopes and R.~Pardini.
\newblock The geography of irregular surfaces.
\newblock {\em Math. Sci. Res. Inst. Publ.}, 59:349–378, 2012.

\bibitem{Par1}
R.~Pardini.
\newblock Abelian covers of algebraic varieties.
\newblock {\em J. Reine Angew. Math.}, 417:191–213, 1991.

\bibitem{Per2}
U.~Persson.
\newblock A family of genus two fibrations.
\newblock {\em Springer Lecture Notes}, 732:496--503, 1979.

\bibitem{Per1}
U.~Persson.
\newblock Chern invariants of surfaces of general $\text{type}$.
\newblock {\em Compositio Mathematica}, 43:3--58, 1981.

\bibitem{PerssonGang}
U.~Persson, C.~Peters, and G.~Xiao.
\newblock Geography of spin surfaces.
\newblock {\em Topology 35}, 4:845–862, 1996.

\end{thebibliography}
\vspace{5mm}

\noindent Vicente Lorenzo \footnote{The author is a doctoral student of the Department of Mathematics and
Center for Mathematical Analysis, Geometry and Dynamical Systems of Instituto Superior T\'{e}cnico,
Universidade de Lisboa and is supported by Fundac\~{a}o para a Ci\^{e}ncia e a Tecnologia (FCT), Portugal through
the program Lisbon Mathematics PhD (LisMath), scholarship  FCT - PD/BD/128421/2017 and
projects UID/MAT/04459/2019 and UIDB/04459/2020.}\\
Center for Mathematical Analysis, Geometry and Dynamical Systems\\
Departamento de Matem\'{a}tica\\
Instituto Superior T\'{e}cnico\\
Universidade de Lisboa\\
Av. Rovisco Pais\\
1049-001 Lisboa\\
Portugal\\
\textit{E-mail address: }{vicente.lorenzo@tecnico.ulisboa.pt}\\
https://orcid.org/0000-0003-2077-6095

\end{document}